\documentclass[11pt, a4paper]{article}
\usepackage{a4wide}
\usepackage{mathrsfs} %mathscr
\usepackage{amssymb} %Box
\usepackage{amsmath} %Operator
\usepackage{amsthm} %proof
\usepackage[hidelinks]{hyperref}

\newtheorem{theorem}{Theorem}
\newtheorem{proposition}[theorem]{Proposition}
\newtheorem{lemma}[theorem]{Lemma}
\newtheorem{corollary}[theorem]{Corollary}

\newcommand{\G}[2]{G_{#1,#2}}
\newcommand{\calG}[2]{{\cal G}_{#1,#2}}
\newcommand{\Exp}{\,\mathbb{E}}
\renewcommand{\Pr}{\,\mathbb{P}}
\newcommand{\eps}{\varepsilon}

\DeclareMathOperator{\Bin}{Bin}
\DeclareMathOperator{\avgdeg}{\overline{\deg}}

\title{Impropriety and resilience of the chromatic number}

\author{
Ross J. Kang
\thanks{Korteweg--de Vries Institute for Mathematics, University of Amsterdam, P.O. Box 94248, 1090 GE Amsterdam, Netherlands.
Email: \protect\href{mailto:r.kang@uva.nl}{\protect\nolinkurl{r.kang@uva.nl}}. Partially supported, at different stages, by Veni (639.031.138) and Vidi (639.032.614) grants of the Netherlands Organisation for Scientific Research (NWO), by the Gravitation Programme NETWORKS (024.002.003) of the Dutch Ministry of Education, Culture and Science (OCW), and by NWO Open Competition grant OCENW.M20.009.}
\and
Colin McDiarmid
\thanks{Department of Statistics, University of Oxford, 24--29 St Giles', Oxford OX1 3LB, United Kingdom.
Email: \protect\href{mailto:cmcd@stats.ox.ac.uk}{\protect\nolinkurl{cmcd@stats.ox.ac.uk}}.}
}

\begin{document}

\maketitle

\begin{abstract}
We show that, for $p$ bounded away from $1$ and $np\to\infty$, the chromatic number of the binomial random graph $G_{n,p}$ is resilient to removal of a subgraph $H$ with $o(n \log(np))$ edges. It follows that the same assertion holds with $H$ being a subgraph of maximum degree $o(\log(np))$; and both these results are best possible up to a constant factor. We also prove corresponding results for the random $r$-regular graph $G_{n,r}$, when $r\to\infty$ and $r/n$ is bounded away from $1$, where $H$ has $o(n\log r)$ edges or maximum degree $o(\log r)$. Our results rely on a natural connection between chromatic resilience to edge removal and improper graph colouring. This work provides counterpoint to the work of Sudakov and Vu (2008) and Alon and Sudakov (2010), which treated chromatic resilience for $G_{n,p}$ subject to the addition, rather than removal, of edges.
\end{abstract}

%%%%%%%%%%%%%%%%%%%%%%%%%%%%%%%%%%%%%%%%%%%%%%%%%%%%%%%%%%%%%%%%%%%%%%%%%%%%%%%%%%%%%%%%%%%%%%%%%%%%%%%%%%%%%%%%%%%%%%%%%%%%%%%%%%%%%%%%%%%%%%%%%%%%%%%%%%%%%%%%%%%%%%%%%%%%%%%%%%%%%%%%%%%%%%%%%%%%%%%%%%%%%%%%%%%%%%%%%%%%%%%%%%%%%%%%%%%%%%%%%%%%%%%%%%%%%%%%%%%%%%%%%%%%%%%%%%%%%%%%%%%%%
\section{Introduction}\label{sec:intro}

Let $\G{n}{p}$ denote the binomial random graph with vertex set $[n] = \{1, \ldots, n \}$ in which the edges are included independently at random with probability $p = p(n) \in (0,1)$ and let $b = 1/(1-p)$.
The chromatic number $\chi(\G{n}{p})$ of $\G{n}{p}$ is a central notion in random graph theory.  The problem of determining this parameter was posed in 1960 by Erd\H{o}s and R\'enyi~\cite{ErRe60}. Their original question (with $np$ a constant) remains open, but there are very strong estimates --- see~\cite{Coj13,CoVi16} for progress.
In the case where $np\to \infty$ and $p$ is bounded away from $1$, it is a classic result due to Bollob\'as~\cite{Bol88}, Matula and Ku\v cera~\cite{MaKu90} and {\L}uczak~\cite{Luc91a} that
\begin{align*}
\chi(\G{n}{p}) = (1+o(1))\frac{n}{2\log_b(np)}
\end{align*}
asymptotically almost surely (a.a.s.), that is, with probability approaching $1$ as $n\to\infty$.
Recall that if $p=o(1)$, then $\log b=(1+o(1))p$, so the right-hand side may be written as $(1+o(1))np/(2\log (np))$.
Throughout, $\log$ denotes the natural logarithm.
For further background into the colouring of random graphs,
consult~\cite{Bol01,JLR00,KaMc15}.

Alon and Sudakov~\cite{AlSu10}, improving upon a result of Sudakov and Vu~\cite{SuVu08}, showed that for dense enough $\G{n}{p}$, specifically $n^{-1/3+\eta}\le p\le 1/2$ for any fixed $\eta > 0$, the following result holds:
a.a.s.~$\chi(\G{n}{p})$
is resilient to the addition of any subset $H$ of
$o(n^2/(\log_b(np))^2)=o((np)^2/(\log n)^2)$
edges on $[n]$;
that is, for any $\eps>0$, there exists $\delta>0$ such that no matter how an adversary chooses a graph $H$ on $[n]$ with at most $\delta n^2/(\log_b(np))^2$ edges, 
\begin{align*}
\chi(\G{n}{p}\cup H) \le (1+\eps)\chi(\G{n}{p}) \text{ a.a.s.}
\end{align*}
This is best possible up to a constant factor, as can be seen by considering $H$ to be the edges of a clique on $(1+\eps)n/(2\log_b(np))$ vertices.
Under the same conditions for $p$, they also proved a local version, for which $H$ has maximum degree
$o(n/(\log_b(np)\log\log n))=o(np/(\log n\log\log n))$;
this is best possible up to a $\log\log n$ factor.

Sudakov and Vu began this line of study within a broader study of resilience of graph properties to removal or addition of edges.
Besides chromatic number, they analysed the resilience in $\G{n}{p}$ of perfect matchings, Hamiltonian cycles and asymmetry, extending work in~\cite{KSV02,KiVu04}, see also~\cite{Sud17} for a review of developments.
In the present work, we are concerned with resilience of the chromatic number to edge {\em removal} rather than addition.
We settle this using a natural relation to improper graph colourings, and show the following.

\begin{theorem}\label{thm:main}
Let $p=p(n)$ be bounded away from $1$ and $np \to\infty$.
The following hold a.a.s.
\begin{enumerate}
\item
$\chi(\G{n}{p})$
is (globally) resilient to removal of any subgraph with $o(n \log(np))$ edges.
\item
$\chi(\G{n}{p})$ is (locally) resilient to removal of any subgraph of maximum degree $o(\log(np))$.
\end{enumerate}
The local assertion is best possible up to a constant factor. The global assertion is best possible.
\end{theorem}
\noindent
Theorem~\ref{thm:main} can be reformulated as follows. For any $\eps>0$, there exists $\delta>0$ such that no matter how an adversary chooses $H \subseteq \G{n}{p}$ with at most $\delta n \log(np)$ edges,
\begin{align*}
\chi(\G{n}{p}-H) \ge (1-\eps)\chi(\G{n}{p}) \text{ a.a.s.:}
\end{align*}
this is the global resilience result, which implies the local one.
On the other hand, for some $\eps,\delta>0$, there exists $H' \subseteq \G{n}{p}$ of maximum degree at most $\delta \log(np)$ such that
\begin{align*}
\chi(\G{n}{p}-H') < (1-\eps)\chi(\G{n}{p}) \text{ a.a.s.;}
\end{align*}
and, for any $\eps > 0$, there exist $\delta > 0$ and $H'' \subseteq \G{n}{p}$ with at most $\eps n\log(np)$ edges such that
\begin{align*}
\chi(\G{n}{p}-H'') < (1-\delta)\chi(\G{n}{p}) \text{ a.a.s.}
\end{align*}

The main message to take from Theorem~\ref{thm:main} is that decreasing $\chi(\G{n}{p})$ is a global task, of removing edges throughout $\G{n}{p}$. In contrast, the results of Alon and Sudakov say that increasing $\chi(\G{n}{p})$ is essentially a local task, of adding a clique, if $p$ does not decay too rapidly.

If the chromatic resilience results for edge addition were extended to the entire range of $p=p(n)$ with $p\to 0$ and $np\to\infty$, then, at least when $np\gg(\log(np))^2\log\log n$, they, together with Theorem~\ref{thm:main}, would imply that, locally, the chromatic number is more resilient to increase than to decrease; whereas, globally, the chromatic number is much more resilient to decrease than to increase if $\G{n}{p}$ is very sparse and the other way around if $\G{n}{p}$ is very dense (with the switch at around $np=\sqrt{n}$, up to logarithmic factors).

We also consider random regular graphs. Their study helped to initiate the broader theory of resilience~\cite{KiVu04}. Let $\G{n}{r}$ denote the random $r$-regular graph, chosen uniformly from the set $\calG{n}{r}$ of all $r$-regular graphs with vertex set $[n]$ (assuming $rn$ is even).
The Kim--Vu sandwich conjecture~\cite{KiVu04} gives a useful connection between the binomial and random regular models. It asserts that, when $r=\omega(\log n)$, one may couple $\G{n}{r}$ with two binomial random graphs $\G{n}{p_-}$ and $\G{n}{p_+}$ such that
\begin{align*}
\G{n}{p_-}\subseteq\G{n}{r}\subseteq\G{n}{p_+}
\end{align*}
a.a.s., where $p_-=(1+o(1))r/n$ and $p_+=(1+o(1))r/n$.
Gao, Isaev and McKay~\cite{GaIsMc25} recently proved it whenever $\min\{r,n-r\}\gg(\log n)^4$.
Two contemporaneous works then treated the remaining regimes. Gao, Isaev and Perez-Gimenez~\cite{GaIsPeGi25} obtained related couplings that settled the remaining very dense regime and gave a weaker form of the conjecture in the remaining sparse regime, while Behague, Il'kovi{\v c} and Montgomery~\cite{BeIlMo25} proved the conjecture in full. We use these sandwich results to prove the following.

\begin{theorem}\label{thm:regular}
Let $r = r(n)$ satisfy $r\to \infty$ and $r\le cn$ for some fixed $c<1$.
The following hold a.a.s.
\begin{enumerate}
\item
$\chi(\G{n}{r})$ is (globally) resilient to removal of any subgraph with $o(n \log r)$ edges.
\item
$\chi(\G{n}{r})$ is (locally) resilient to removal of any subgraph of maximum degree $o(\log r)$.
\end{enumerate}
The local assertion is best possible up to a constant factor. The global assertion is best possible in the sense described after Theorem~\ref{thm:main}.
\end{theorem}

The structure of the paper is as follows.
We show links between chromatic resilience and several improper colouring parameters.
We use these to prove Theorem~\ref{thm:main} in Section~\ref{sec:main} and Theorem~\ref{thm:regular} in Section~\ref{sec:regular}.
For the regular result, we combine a direct pairing argument with the sandwich result mentioned above. We also give sharper fixed-density forms (Theorem~\ref{thm:sharp,dense} and Corollary~\ref{cor:sharp,regular,dense}) in Section~\ref{sec:further}.

\section{Impropriety and resilience}\label{sec:impropriety}

As mentioned above, we can relate the chromatic number's (removal) resilience qualities with improper colouring parameters, one of which was already studied in the 1960s~\cite{Lov66}.
The {\it (local) $t$-improper chromatic number} $\chi^t(G)$ of a graph $G$ is the smallest number of colours needed in a {\it (local) $t$-improper colouring} --- a colouring of the vertices of $G$ so that each subgraph induced by some colour class has maximum degree at most $t$.
See~\cite{Woo18} for a survey of this type of colouring parameter, especially with respect to sparse graph classes.
The {\it global $t$-improper chromatic number} $\hat{\chi}^t(G)$ of a graph $G$ on $n$ vertices is the smallest number of colours needed in a {\it global $t$-improper colouring} --- a colouring of the vertices of $G$ so that the total number of edges in the subgraphs induced by colour classes is at most $nt/2$.
Equivalently, the spanning subgraph of monochromatic edges has average degree at most $t$. We write $\avgdeg(G)$ for the average degree of $G$. The global $t$-improper chromatic number appears to have received little attention in graph theory. It was suggested at the end of a paper from 1997~\cite{CGJ97}. An equivalent formulation in terms of chromatic edge-deletion blockers has been recently considered from a computational perspective~\cite{LuMa24}.

Here and below the parameter $t$ is a non-negative integer for local impropriety, while it may be any non-negative real number for global impropriety.

Let us make a few remarks on the basic behaviour of these parameters.
Both of them are non-increasing in $t$ and include the chromatic number as a special case by taking $t=0$.  Clearly $\hat{\chi}^t(G)\le\chi^{\lfloor t\rfloor}(G)$ for every real $t\ge0$.  Furthermore, it follows from a result of Lov\'asz~\cite{Lov66} that $\chi^t(G) \le \lceil(\Delta(G)+1)/(t+1)\rceil$ for every non-negative integer $t$, where $\Delta(G)$ denotes the maximum degree of $G$.  Summarising, we have that for any graph $G$ and non-negative real $t$
\begin{align}\label{eqn:basic}
\hat{\chi}^t(G) \le \chi^{\lfloor t\rfloor}(G) \le \min\left\{ \left\lceil \frac{\Delta(G)+1}{\lfloor t\rfloor+1} \right\rceil, \chi(G) \right\}.
\end{align}

We now present two elementary correspondences between local/global resilience of the chromatic number and local/global improper colouring.

\begin{proposition}\label{prop:localcorr}
Let $G$ be a graph, $t$ a non-negative integer, and $k$ a positive integer. The following statements are equivalent:
\begin{enumerate}
\item
$\chi^t(G) \le k$.
\item
There exists a subgraph $H$ of $G$ with $\Delta(H)\le t$ such that $\chi(G-H) \le k$.
\end{enumerate}
\end{proposition}

\begin{proposition}\label{prop:globalcorr}
Let $G$ be a graph on $n$ vertices, $t$ a non-negative real number, and $k$ a positive integer. The following statements are equivalent:
\begin{enumerate}
\item
$\hat{\chi}^t(G) \le k$.
\item
There exists a subgraph $H$ of $G$ with at most $nt/2$ edges such that $\chi(G-H) \le k$.
\end{enumerate}
\end{proposition}

Both correspondences follow by deleting the monochromatic edges of an improper colouring and, conversely, using a proper colouring of $G-H$.

Note that the correspondences of Propositions~\ref{prop:localcorr} and~\ref{prop:globalcorr} apply to all graphs and allow us to deduce chromatic resilience results from prior work on improper chromatic numbers.

Given the above correspondences, it can be helpful to relate chromatic resilience to two other parameters that are also closely related to improper colouring, particularly for random graphs, cf.~\cite{BoCl95,KaMc10,FKM10,FKM14}.
The {\em (local) $t$-stability number~$\alpha^t(G)$} of a graph $G$ is the number of vertices in a largest {\em (local) $t$-stable set}, a subset of the vertices whose induced subgraph has maximum degree at most $t$.
The {\em global $t$-stability number~$\hat{\alpha}^t(G)$} of a graph $G$ is the number of vertices in a largest {\em global $t$-stable set}, a subset of the vertices whose induced subgraph has average degree at most $t$. (We called global $t$-stability {\em $t$-sparsity} in~\cite{FKM14}.)

Both of these parameters are non-decreasing in $t$ and include the stability number as a special case by taking $t=0$.
For every non-negative real $t$, we have $\alpha^{\lfloor t\rfloor}(G) \le \hat{\alpha}^t(G)$ and moreover
\begin{align}\label{eqn:basic2}
\chi(G) \ge \chi^{\lfloor t\rfloor}(G) \ge \frac{|V(G)|}{\alpha^{\lfloor t\rfloor}(G)} \ge \frac{|V(G)|}{\hat{\alpha}^t(G)}.
\end{align}
Thus we may bound $\chi^{\lfloor t\rfloor}(G)$ from below by bounding $\hat{\alpha}^t(G)$ or $\alpha^{\lfloor t\rfloor}(G)$ from above.

As we next show, we can also bound $\hat{\chi}^t(G)$ from below by bounding the global stability of $G$ from above, although the connection holds in a weaker sense than above in that we need to bound $\hat{\alpha}^{t'}(G)$ for some $t'>t$. (Note that in this case $\hat{\alpha}^{t'}(G) \ge \hat{\alpha}^{t}(G)$ always.)

\begin{proposition}
\label{prop:stability}
Let $G$ be a graph on $[n]$ and $\eps > 0$.
Suppose that there exist $t\ge0$ and $\alpha>0$ such that $\hat{\alpha}^{2t/\eps}(G) < (1+\eps/2)\alpha$.
Then $\hat{\chi}^t(G) > \lfloor (1-\eps)n/\alpha \rfloor$.
\end{proposition}

\begin{proof}
Let $k = \lfloor (1-\eps)n/\alpha \rfloor$. If $k\le0$, the result is immediate. Otherwise, suppose for a contradiction that $\hat{\chi}^t(G) \le k$. Consider any global $t$-improper colouring $P$ of the graph using $k$ colours and let $A\subseteq [n]$ be the union of those colour classes of $P$ having fewer than $(1+\eps/2)\alpha$ vertices. We then have
\begin{align*}
|A| \le k(1+\eps/2)\alpha \le (1-\eps)(1+\eps/2)n < (1-\eps/2)n.
\end{align*}
It follows that the number of vertices in colour classes of $P$ having at least $(1+\eps/2)\alpha$ vertices is more than $\eps n/2$.
Since $\hat{\alpha}^{2t/\eps}(G) < (1+\eps/2)\alpha$, each such colour class of $P$ induces a subgraph of average degree greater than $2t/\eps$.
Summing over all these vertices, we see that the number of monochromatic edges under $P$ is more than
$(\eps n/2)(2t/\eps)/2 = nt/2$,
contradicting that $P$ is a global $t$-improper colouring.
\end{proof}

\section{Proof of Theorem~\ref{thm:main}}\label{sec:main}

We now apply Propositions~\ref{prop:localcorr} and~\ref{prop:globalcorr} to the setting of random graphs.  Given these correspondences, it suffices for the proof of our first main result to establish appropriate bounds on the local and global $t$-improper chromatic numbers.

\begin{lemma}\label{lem:lower}
Suppose $p = p(n)$ is bounded away from $1$ and $\eps > 0$.  There exist constants $d_0$ and $\tau > 0$ such that, if $np \ge d_0$ and $0\le t \le \tau \log (np)$, then $\hat{\chi}^t(\G{n}{p}) \ge n/(2(1+\eps)\log_b (np))$ a.a.s.
\end{lemma}

\begin{lemma}\label{lem:upper}
Suppose $p = p(n)$ is bounded away from $1$, $np\to\infty$ and $\eps > 0$.  If $t=t(n)$ is a positive integer, then $\chi^t(\G{n}{p}) \le (1+\eps)np/t + np/(\log(np))^2$ a.a.s.
\end{lemma}

With these lemmas in hand, let us see how to prove the main theorem.

\begin{proof}[Proof of Theorem~\ref{thm:main}]
We first prove global resilience --- from this, the local assertion follows since $\hat{\chi}^t(\G{n}{p})\le\chi^t(\G{n}{p})$.  Let $t=t(n)=o(\log(np))$.  Proposition~\ref{prop:globalcorr} and Lemma~\ref{lem:lower} give, a.a.s.~and simultaneously for every subgraph $H$ with at most $nt/2$ edges,
\begin{align*}
\chi(\G{n}{p}-H)\ge(1+o(1))\frac{n}{2\log_b(np)}.
\end{align*}
Together with the usual upper bound on $\chi(\G{n}{p})$, this proves the first two assertions.

We next prove that the local assertion is best possible up to a constant factor.  Take
\begin{align*}
t=\left\lceil4p\log_b(np)\right\rceil.
\end{align*}
Uniformly for $p$ bounded away from $1$, we have $t=\Theta(\log(np))$.  Moreover,
\begin{align*}
\frac{np}{(\log(np))^2}=o\left(\frac{n}{\log_b(np)}\right).
\end{align*}
Lemma~\ref{lem:upper}, with $\eps=1/4$, therefore gives
\begin{align*}
\chi^t(\G{n}{p})<\frac{n}{3\log_b(np)}
\end{align*}
for all sufficiently large $np$.  The required subgraph now follows from Proposition~\ref{prop:localcorr}.

It remains to prove that the global assertion is best possible.  Fix $0<\eps<1/2$ and set $\delta=\eps/3$.  We use the standard a.a.s.~estimates
\begin{align*}
\chi(\G{n}{p})&=(1+o(1))\frac{n}{2\log_b(np)},&
\alpha(\G{n}{p})&=(1+o(1))2\log_b(np)
\end{align*}
simultaneously; see, for example,~\cite{Bol01,JLR00}.  We shall also use the following density property.  If $\xi=(np)^{-1/4}$, then a Chernoff bound followed by a union bound over the at most $2^n$ sets shows that, for any fixed $0<a_1<a_2\le1$, a.a.s.~simultaneously for every $S\subseteq[n]$ with $a_1n\le|S|\le a_2n$,
\begin{align*}
e(S,[n]\setminus S)+e(S)
=(1+O(\xi))p\left(|S|(n-|S|)+\binom{|S|}{2}\right).
\end{align*}
Indeed, the exponent in the Chernoff bound has order $\xi^2n^2p=n\sqrt{np}$.

Fix an optimal proper colouring with $q=\chi(\G{n}{p})$ colours, take $s=\lfloor\delta q\rfloor$, and let $U$ be the union of the $s$ largest colour classes. Writing the sizes of the colour classes in non-increasing order as $c_1,\ldots,c_q$, we have
\begin{align*}
|U|=\sum_{i=1}^s c_i
\ge \frac{s}{q}\sum_{i=1}^q c_i
=\frac{sn}{q}
=(\delta+o(1))n.
\end{align*}
On the other hand, each colour class has size at most $\alpha(\G{n}{p})$, and hence
\begin{align*}
|U|\le s\alpha(\G{n}{p})=(\delta+o(1))n.
\end{align*}
Consequently, $|U|=(\delta+o(1))n$.
Retain the other $q-s$ colour classes and order the vertices of $U$.  Recolour them in this order, each time choosing a retained colour that minimises the number of neighbours already given that colour.  Every resulting monochromatic edge within $U$ is then counted when its later endpoint is recoloured.  The total number of monochromatic edges is consequently at most
\begin{align*}
\frac{e(U,[n]\setminus U)+e(U)}{q-s}
&=(1+o(1))\frac{2\delta-\delta^2}{1-\delta}\,np\log_b(np)\\
&\le\eps n\log(np),
\end{align*}
where we used $p/\log b\le1$.  Removing these monochromatic edges leaves a graph coloured with $(1-\delta+o(1))\chi(\G{n}{p})$ colours, which is smaller than $(1-\delta/2)\chi(\G{n}{p})$ for all sufficiently large $n$. This proves the assertion.
\end{proof}

Before giving the proofs of Lemmas~\ref{lem:lower} and~\ref{lem:upper}, we introduce some notation and a lemma from the theory of large deviations.
For further background, consult~\cite{DeZe98}.  Given $p\in(0,1)$, let
\begin{align}
\label{eqn:Lambda*}
\Lambda^*(x) = \left\{ \begin{array}{ll}
\displaystyle x \log \frac{x}{p} + (1 - x) \log \frac{1 - x}{1-p} & \mbox{for $x\in [0, 1]$}\\
\infty & \mbox{otherwise}
\end{array} \right.
\end{align}
(where $\Lambda^*(0) = \log b$ and $\Lambda^*(1) = -\log p$).
This is the Fenchel-Legendre transform of the logarithmic moment generating function associated with the Bernoulli distribution with probability $p$ (cf.~\cite[Exercise~2.2.23(b)]{DeZe98}).  It is easy to see that $\Lambda^*(x)$ has a global minimum of $0$ at $x = p$, is strictly decreasing on $[0, p)$ and strictly increasing on $(p,1]$.  
Let $\Bin(\cdot,\cdot)$ denote the binomial distribution.
We use the following estimate from~\cite[Lemma~4]{KaMc10}.

\begin{lemma}\label{lem:binomial}
There is a constant $\delta>0$ such that the following holds.  Let $p\in(0,1)$, let $N$ be a positive integer,
and let $X \in \Bin(N,p)$.  Then, for each integer $1\le r\le N-1$ such that $r\le Np$,
\[
\delta \cdot\max\left\{r^{-1/2},(N-r)^{-1/2}\right\} \cdot \exp(-N\Lambda^*( r/N )) \le \Pr( X \le r) \le \exp(-N\Lambda^*( r/N )).
\]
\end{lemma}

\begin{lemma}\label{lem:rate}
For every $c_0,\zeta>0$, there is some $a_0>0$ such that, for $0<p\le 1-c_0$ and $0\le a\le a_0$,
\begin{align*}
\Lambda^*(ap)\ge (1-\zeta)\log b.
\end{align*}
\end{lemma}

\begin{proof}
We first consider $p$ close to zero.  Uniformly for $a$ in any fixed interval $[0,A]$, with $a\log a$ interpreted as zero when $a=0$, we have
\begin{align*}
\Lambda^*(ap)
&=ap\log a+(1-ap)\log\frac{1-ap}{1-p}\\
&=p(1-a+a\log a)+O(p^2).
\end{align*}
Since $\log b=p+O(p^2)$ and $1-a+a\log a\to1$ as $a\to0$, the required inequality follows for all sufficiently small $p$ and $a$.

For $p$ in a compact interval bounded away from zero, the assertion follows from uniform continuity at $a=0$, since $\Lambda^*(0)=\log b$.  Combining the two ranges gives a choice of $a_0$ that is uniform over $0<p\le1-c_0$.
\end{proof}

Before continuing, we remark that an alternative way to establish Lemma~\ref{lem:lower} is to apply Lemma~\ref{lem:binomial} to obtain an upper bound on the global $t$-stability number of $\G{n}{p}$ and then apply Proposition~\ref{prop:stability}; however, this gives a smaller admissible value of $\tau$. The method we demonstrate here is sharp (with respect to $\tau$) for fixed $p$, as we show in Section~\ref{sec:further}. 

\begin{proof}[Proof of Lemma~\ref{lem:lower}]
Choose $\zeta>0$ so that $(1+\eps)(1-\zeta)>1$.  Let $k=2(1+\eps)\log_b(np)$ and $j=\lfloor n/k\rfloor$.  Consider any colouring of $[n]$ using at most $j$ colours, and let $M$ be the number of pairs of vertices contained in its colour classes.  By convexity,
\begin{align*}
M\ge\frac12\left(\frac{n^2}{j}-n\right)\ge\frac{n(k-1)}2.
\end{align*}

Let $N_0=\lfloor n(k-1)/2\rfloor$ and $T=\lfloor nt/2\rfloor$.  For $t\le\tau\log(np)$, we have
\begin{align*}
\frac{T}{pN_0}
\le (1+o(1))\frac{\tau}{2(1+\eps)}\frac{\log b}{p}.
\end{align*}
Since $\log b/p$ is bounded uniformly for $p$ bounded away from $1$, we may choose $\tau$ so that, for all sufficiently large $n$, we have $T\le pN_0<N_0$ and Lemma~\ref{lem:rate} yields
\begin{align*}
\Lambda^*\left(\frac{T}{N_0}\right)\ge(1-\zeta)\log b.
\end{align*}
The number of monochromatic edges in the colouring stochastically dominates $\Bin(N_0,p)$.  For $T\ge1$, Lemma~\ref{lem:binomial} therefore gives
\begin{align*}
\Pr\big(\text{the colouring is globally $t$-improper}\big)
\le\exp\left(-N_0\Lambda^*\left(\frac{T}{N_0}\right)\right).
\end{align*}
If $T=0$, the colouring is globally $t$-improper only if none of its $M$ candidate monochromatic pairs forms an edge. Since $M\ge N_0$, the probability of this event is at most
\begin{align*}
(1-p)^{N_0}
=\exp\bigl(-N_0\Lambda^*(0)\bigr),
\end{align*}
so the same estimate remains valid.

Moreover,
\begin{align*}
j\le\frac{n\log b}{2(1+\eps)\log(np)},
\end{align*}
and hence $\log j\le\log(np)-\log\log(np)+O(1)$.  A union bound over the at most $j^n$ labelled colourings now gives a probability at most
\begin{align*}
\exp\left(-\big((1+\eps)(1-\zeta)-1+o(1)\big)n\log(np)
-n\log\log(np)+O(n)\right)=o(1).
\end{align*}
Thus a.a.s.~there is no global $t$-improper colouring using at most $j$ colours, as required.
\end{proof}

\begin{proof}[Proof of Lemma~\ref{lem:upper}]
Let $\eta=\min\{\eps,1\}/100$.  Choose $C$ sufficiently large.  If $np\ge C\log n$, a Chernoff bound and a union bound over the vertices give $\Delta(\G{n}{p})\le(1+\eta)np$ a.a.s.  The result then follows directly from~\eqref{eqn:basic} for all sufficiently large $n$.

It remains to consider $np<C\log n$.  Let
\begin{align*}
W&=\{v:\deg(v)\ge(1+4\eta)np\},\\
k&=\lceil\eta np\rceil,\qquad
W'=\{v\in W:\deg_{\G{n}{p}[W]}(v)\ge k\}.
\end{align*}
We shall show that the graph induced by $W'$ has degeneracy less than
\begin{align*}
q=\left\lfloor\frac{np}{4(\log(np))^2}\right\rfloor.
\end{align*}

We first bound the size of $W'$.  If $v\in W'$, it has some set $U\subseteq N(v)\cap W$ of size $k$.  Condition on the edges from $v$ to $U$, and for each $u\in U$ let
\begin{align*}
X_u=e\bigl(\{u\},[n]\setminus(\{v\}\cup U)\bigr).
\end{align*}
The edge indicators defining the variables $X_u$, $u\in U$, are disjoint, and hence these variables are independent.  Since $u\in W$, we must have $X_u\ge(1+2\eta)np$ for all sufficiently large $np$.  It follows from a Chernoff bound, cf.~\cite[Theorem~2.1]{JLR00}, that for some constants $c,c'>0$,
\begin{align*}
\Pr(v\in W')
&\le\binom{n-1}{k}p^k\exp(-c\eta^2npk)
\le\exp(-c'(np)^2).
\end{align*}
Let $w_0=n\exp(-c'(np)^2/2)$.  If $w_0\ge1$, then Markov's inequality gives $|W'|\le w_0$ a.a.s.; if $w_0<1$, then $\Pr(W'\ne\emptyset)\le\Exp|W'|=o(1)$.  Thus we may assume that $|W'|\le w_0$.

We next rule out all dense subsets of this size simultaneously.  For $q+1\le s\le w_0$, let $m=\lceil qs/2\rceil$.  Since $\log(n/s)\ge c'(np)^2/2$, a union bound gives
\begin{align*}
\Pr\big(\text{some $S$, $|S|=s$, has $e(S)\ge m$}\big)
&\le\binom ns\binom{\binom s2}{m}p^m\\
&\le\left(\frac{en}{s}\right)^s
     \left(\frac{esp}{q}\right)^{qs/2}
\le\exp\left(-\frac{qs}{3}\log\frac ns\right)
\end{align*}
for all sufficiently large $np$.  Summing this estimate over $s$ gives $o(1)$.  Every graph of minimum degree at least $q$ contains at least $q+1$ vertices and at least $qs/2$ edges, so $\G{n}{p}[W']$ has degeneracy at most $q-1$ a.a.s.

Finally, use disjoint sets of colours on $[n]\setminus W$, $W\setminus W'$ and $W'$.  By~\eqref{eqn:basic} and the degeneracy bound,
\begin{align*}
\chi^t(\G{n}{p})
&\le\left\lceil\frac{(1+4\eta)np+1}{t+1}\right\rceil
  +\left\lceil\frac{k}{t+1}\right\rceil+q\\
&\le(1+6\eta)\frac{np}{t}+2+q
\le(1+\eps)\frac{np}{t}+\frac{np}{(\log(np))^2},
\end{align*}
as required.
\end{proof}

\section{Proof of Theorem~\ref{thm:regular}}\label{sec:regular}

The results of Frieze and \L uczak~\cite{FrLu92} and Cooper, Frieze, Reed and Riordan~\cite{CFRR02} give the ordinary chromatic and stability number estimates needed below when $r$ is at most a constant multiple of $\log n$.  Recently, Behague, Il'kovi{\v c} and Montgomery~\cite{BeIlMo25} proved the Kim--Vu sandwich conjecture in full.  Their quantitative result transfers the corresponding binomial estimates above that range.  Thus, whenever $r\to\infty$ and $r\le cn$ for some fixed $c<1$, we may use
\begin{align}\label{eqn:regular-baseline}
\chi(\G{n}{r})&=(1+o(1))\frac{n}{2\log_b r},&
\alpha(\G{n}{r})&=(1+o(1))2\log_b r
\end{align}
a.a.s., where $b=n/(n-r)$.

Recall~\eqref{eqn:basic2}.  For the lower range, we shall need the following upper bound on the global $t$-stability number.

\begin{lemma}\label{lem:reg,local}
Fix $C,\eps>0$.  Suppose $r=r(n)$ is a positive integer, $rn$ is even, $r\to\infty$ and $r\le C\log n$.  There is some $\tau>0$ such that, if $0\le t\le\tau\log r$, then
\begin{align*}
\hat{\alpha}^t(\G{n}{r})<(1+\eps)2\log_b r
\end{align*}
a.a.s., where $b=n/(n-r)$.
\end{lemma}

Let us first see how this lemma is used to prove Theorem~\ref{thm:regular}.

\begin{proof}[Proof of Theorem~\ref{thm:regular}]
Fix $\rho>0$, and let $t=t(n)=o(\log r)$.  Choose $\eta>0$ sufficiently small in terms of $\rho$ and the constant $c$ in the theorem.  By~\cite[Theorem~1.1]{BeIlMo25}, there is some $C=C(\eta)$ such that, whenever $r\ge C\log n$, there are probabilities $p_-\le p_+$ and a coupling for which
\begin{align}\label{eqn:regular-sandwich}
\G{n}{p_-}\subseteq\G{n}{r}\subseteq\G{n}{p_+}
\end{align}
a.a.s., where
\begin{align*}
(1-\eta)\frac rn\le p_-\le p_+\le(1+\eta)\frac rn.
\end{align*}
We choose $\eta$ small enough that $(1+\eta)c<1$.

Suppose first that $r<C\log n$.  Apply Lemma~\ref{lem:reg,local} with $\eps=\rho/6$ and with the parameter $t$ in that lemma equal to $6t/\rho$.  Since $t=o(\log r)$, this gives
\begin{align*}
\hat{\alpha}^{6t/\rho}(\G{n}{r})<(1+\rho/6)2\log_b r
\end{align*}
a.a.s.  Proposition~\ref{prop:stability}, applied with $\eps=\rho/3$ and $\alpha=2\log_b r$, now yields
\begin{align*}
\hat{\chi}^t(\G{n}{r})>\left\lfloor\frac{(1-\rho/3)n}{2\log_b r}\right\rfloor.
\end{align*}

We next suppose that $r\ge C\log n$.  The parameters $\chi^t$ and $\hat{\chi}^t$ are non-decreasing under the addition of edges, so~\eqref{eqn:regular-sandwich} and Lemma~\ref{lem:lower} give the required lower bound from $\G{n}{p_-}$.  For completeness, if $b_\pm=1/(1-p_\pm)$, then uniformly for $r\le cn$,
\begin{align*}
\frac{n}{2\log_{b_\pm}(np_\pm)}
=\big(1+O_c(\eta)+o(1)\big)\frac{n}{2\log_b r}.
\end{align*}
The ordinary upper bound follows in the same way from
$\chi(\G{n}{r})\le\chi(\G{n}{p_+})$.  Since $\rho$ was arbitrary, Propositions~\ref{prop:localcorr} and~\ref{prop:globalcorr} prove both resilience assertions.

It remains to prove the two sharpness statements.  By~\eqref{eqn:basic} and~\eqref{eqn:regular-baseline}, with $t=\lceil8\log r\rceil$,
\begin{align*}
\frac{\chi^t(\G{n}{r})}{\chi(\G{n}{r})}
\le(1+o(1))\frac{r}{4n\log b}\le\frac14+o(1).
\end{align*}
This proves the local assertion by Proposition~\ref{prop:localcorr}.

For the global assertion, fix an optimal proper colouring with $q=\chi(\G{n}{r})$ colours, discard $s=\lfloor\delta q\rfloor$ colour classes, and let $U$ be their union.  Recolour the vertices of $U$ sequentially with the remaining $q-s$ colours, each time minimising the number of already coloured neighbours receiving the chosen colour.  By regularity and~\eqref{eqn:regular-baseline}, the number of resulting monochromatic edges is at most
\begin{align*}
\frac{r|U|}{q-s}
&\le(1+o(1))\frac{\delta}{1-\delta}r\alpha(\G{n}{r})\\
&=(1+o(1))\frac{2\delta}{1-\delta}
  \frac{r/n}{-\log(1-r/n)}\,n\log r\\
&\le(1+o(1))\frac{2\delta}{1-\delta}\,n\log r.
\end{align*}
Given $\eps>0$, choose $\delta>0$ so that $2\delta/(1-\delta)<\eps$. Then, a.a.s., the number of monochromatic edges is at most $\eps n\log r$. Removing these edges leaves a graph coloured with $(1-\delta+o(1))q$ colours, fewer than $(1-\delta/2)q$ for all sufficiently large $n$, as required.
\end{proof}

We next prove Lemma~\ref{lem:reg,local}. Part of this result appears in a non-refereed proceedings contribution by Bollob\'as and Clark~\cite{BoCl95}. Since that work is not easily accessible, we include a full proof.

\begin{proof}[Proof of Lemma~\ref{lem:reg,local}]
Let $k=\lceil2(1+\eps)\log_b r\rceil$.  It is enough to show that a.a.s.~there is no global $t$-stable set of exactly $k$ vertices.  Indeed, if a set $T$ of $K\ge k$ vertices spans at most $tK/2$ edges, then a uniform random $k$-subset $S$ of $T$ satisfies
\begin{align*}
\Exp e(S)
=e(T)\frac{k(k-1)}{K(K-1)}
\le\frac{tk}{2}.
\end{align*}
Thus some such $S$ is globally $t$-stable.

We use the pairing model.  A pairing is a perfect matching of $N=rn$ points partitioned into $n$ cells of $r$ points each, and it induces the usual random regular multigraph.  Fix $k$ cells, write $M=rk$, and let $Z$ denote the number of pairs having both points in these cells.  Writing $P_s=\Pr(Z=s)$, a direct count gives, for
\begin{align*}
\max\{0,M-N/2\}\le s\le\lfloor M/2\rfloor,
\end{align*}
\begin{align*}
P_s
=\binom{M}{2s}(2s-1)!!
  \frac{(N-M)!}{(N-2M+2s)!}
\frac{(N-2M+2s-1)!!}{(N-1)!!}.
\end{align*}
We have $P_s=0$ outside this range. Here and below we use $(-1)!!=1$, corresponding to the unique pairing of the empty set.  Cancelling consecutive terms, we obtain
\begin{align}\label{eqn:pairing-ratio}
\frac{P_{s+1}}{P_s}
=\frac{(M-2s)(M-2s-1)}
       {2(s+1)(N-2M+2s+2)}.
\end{align}

Since $\log b\ge r/n$, we have
\begin{align*}
\frac{M}{N}=\frac{k}{n}\le\frac{2(1+\eps)\log r}{r}+\frac1n=o(1).
\end{align*}
In particular, $M<N/2$ for all sufficiently large $n$. Let
\begin{align*}
\lambda=\frac{M^2}{2(N-2M)}.
\end{align*}
Iterating~\eqref{eqn:pairing-ratio} gives $P_s\le P_0\lambda^s/s!$.  Exposing in turn the partners of the $M$ points also gives
\begin{align*}
P_0
&=\prod_{h=0}^{M-1}\frac{N-M-h}{N-2h-1}\\
&=\prod_{h=0}^{M-1}\left(1-\frac{M-h-1}{N-2h-1}\right)
\le\exp\left(-\frac{M(M-1)}{2N}\right).
\end{align*}

Let $q=\lfloor tk/2\rfloor$. Since $M<N/2$, the support of $Z$ is $\{0,\ldots,\lfloor M/2\rfloor\}$. Moreover,
\begin{align*}
\frac{q}{M/2}\le\frac{t}{r}=o(1),
\end{align*}
so $q$ lies in this support for all sufficiently large $n$. Uniformly for $r\to\infty$ and $r\le C\log n$,
\begin{align*}
k=(1+o(1))\frac{2(1+\eps)n\log r}{r},
\end{align*}
and hence
\begin{align*}
\frac{q}{\lambda}
\le\frac{t(N-2M)}{rM}
\le\frac{tn}{rk}
\le(1+o(1))\frac{t}{2(1+\eps)\log r}
\le(1+o(1))\frac{\tau}{2(1+\eps)}.
\end{align*}
Thus we may choose $\tau$ so that $q\ge1$ implies $q\le\lambda$. We use the bound on $P_0$ directly for $q=0$; otherwise,
\begin{align*}
\Pr(Z\le q)
\le P_0(q+1)\frac{\lambda^q}{q!}
\le P_0(q+1)\left(\frac{e\lambda}{q}\right)^q.
\end{align*}

Let $Y$ count the $k$-sets of cells that contain at most $q$ internal pairs. For a simple pairing, the cells are the vertices of the resulting graph, and the internal pairs are exactly the edges induced by the corresponding vertex set. Consequently, $Y$ counts precisely the global $t$-stable $k$-sets. Let $A=n(\log r)^2/r$. Uniformly for $r\to\infty$ and $r\le C\log n$, we have
\begin{align*}
\log b&=(1+o(1))\frac rn,\\
\log\binom nk
&\le\frac{2(1+\eps)n\log r}{r}\big(\log r-\log\log r+O(1)\big),\\
\frac{M(M-1)}{2N}
&=\big(2(1+\eps)^2+o(1)\big)\frac{n(\log r)^2}{r}.
\end{align*}
Suppose first that $q\ge1$. Then $q=\lfloor tk/2\rfloor\ge tk/4$, and hence
\begin{align*}
q&\le(1+\eps+o(1))\frac{tn\log r}{r},\\
\frac{\lambda}{q}
&\le\frac{2M^2}{(N-2M)tk}
\le C_\eps\frac{\log r}{t}.
\end{align*}
Using $\log(q+1)\le q$ and enlarging $C_\eps$ if necessary, we obtain
\begin{align*}
\log(q+1)+q\log(e\lambda/q)
&\le q\log(e^2\lambda/q)\\
&\le(1+\eps+o(1))\frac{t}{\log r}\log\frac{C_\eps\log r}{t}\,A.
\end{align*}
Together with $\Exp Y\le\binom nk\Pr(Z\le q)$, the estimates above give
\begin{align*}
\log\Exp Y
\le\left(-2\eps(1+\eps)+(1+\eps+o(1))\frac{t}{\log r}\log\frac{C_\eps\log r}{t}+o(1)\right)A.
\end{align*}
For $q=0$, the direct bound $\Pr(Z\le q)=P_0$ instead gives
\begin{align*}
\log\Exp Y\le\big(-2\eps(1+\eps)+o(1)\big)A.
\end{align*}
For $0<t\le\tau\log r$, the term $(t/\log r)\log(C_\eps\log r/t)$ tends to zero uniformly as $\tau\downarrow0$. Since $-2\eps(1+\eps)<0$, we may choose $\tau>0$, depending only on $C$ and the fixed $\eps$, so that, for some $c_\eps>0$, both cases give
\begin{align*}
\Exp Y\le\exp(-c_\eps A).
\end{align*}

Finally, let $\mathcal S$ be the event that the pairing is simple. By Cooper, Frieze, Reed and Riordan~\cite[Lemma~2.2]{CFRR02},
\begin{align*}
\Pr(\mathcal S)\ge\exp(-2r^2)
\end{align*}
provided $100\le r\le n/1000$. Since $r\to\infty$ and $r\le C\log n$, these inequalities hold for all sufficiently large $n$. Moreover,
\begin{align*}
\frac{A}{r^2}
=\frac{n(\log r)^2}{r^3}
\ge\frac{n}{C^3(\log n)^3}
\longrightarrow\infty.
\end{align*}
By Markov's inequality,
\begin{align*}
\Pr(Y>0\mid\mathcal S)
\le\frac{\Pr(Y>0)}{\Pr(\mathcal S)}
\le\frac{\Exp Y}{\Pr(\mathcal S)}
\le\exp(-c_\eps A+2r^2)
=o(1).
\end{align*}
Taking $t=\tau\log r$, write $q_\tau=\lfloor\tau k\log r/2\rfloor$, and let $Y_\tau$ count the $k$-sets of cells that contain at most $q_\tau$ internal pairs. The preceding bound gives
\begin{align*}
\Pr(Y_\tau>0\mid\mathcal S)=o(1).
\end{align*}
If $Y_\tau=0$, then no $k$-set is globally $t$-stable for any real $0\le t\le\tau\log r$. Indeed, such a set would satisfy
\begin{align*}
e(S)\le tk/2\le\tau k\log r/2,
\end{align*}
and hence, since $e(S)$ is an integer, $e(S)\le q_\tau$, contrary to $Y_\tau=0$. Every simple $r$-regular graph arises from exactly $(r!)^n$ pairings, so the pairing conditioned on $\mathcal S$ induces the uniform random graph $\G{n}{r}$, and the result follows.
\end{proof}

\section{Further resilience}\label{sec:further}

One might wonder how resilient $\chi(\G{n}{p})$ is to reduction by a specified non-trivial constant factor. That is, given $\gamma \in(0,1]$, what degree or density conditions on $H\subseteq\G{n}{p}$ ensure that
\begin{align*}
\chi(\G{n}{p}-H) \ge (1+o(1))\gamma \,\chi(\G{n}{p}) \text{ a.a.s.?}
\end{align*}
For fixed $\gamma\in(0,1)$, the arguments for Theorem~\ref{thm:main} give the correct orders $\log(np)$ for $\Delta(H)$ and $n\log(np)$ for $e(H)$, but do not settle the sharp constants.

Our main result in earlier work~\cite{KaMc10}, together with a direct first moment argument for global impropriety, answers this question when $p$ is fixed.  Write $b=1/(1-p)$.  We first define the constant that occurs in the answer.  With $\Lambda^*$ as in~\eqref{eqn:Lambda*}, for $p\in(0,1)$ and $\tau\ge0$, let $\kappa_p(\tau)>\tau/p$ be the unique value such that
\begin{align*}
\frac{\kappa}{2} \Lambda^*\left( \frac{\tau}{\kappa} \right)
\begin{cases}
< 1 & \text{ if }\tau/p < \kappa < \kappa_p(\tau) \\
= 1 & \text{ if }\kappa = \kappa_p(\tau) \\
> 1 & \text{ if }\kappa > \kappa_p(\tau)
\end{cases}
.
\end{align*}
The function $\kappa_p$ is continuous and strictly increasing, with $\kappa_p(0)=2/\log b$.  We define $\tau=\tau(p,\gamma)$ by
\begin{align}\label{eqn:tau-gamma}
\frac{\kappa_p(0)}{\kappa_p(\tau)}=\gamma.
\end{align}

\begin{theorem}\label{thm:sharp,dense}
Fix $p\in(0,1)$ and $\gamma\in(0,1]$, and let $\tau=\tau(p,\gamma)$ be given by~\eqref{eqn:tau-gamma}.  For every $\eps>0$, there is some $\delta>0$ such that the following hold a.a.s.
\begin{enumerate}
\item Every subgraph $H\subseteq\G{n}{p}$ with at most $(\tau/2+\delta)n\log(np)$ edges satisfies
\begin{align*}
\chi(\G{n}{p}-H)\ge(1-\eps)\gamma\chi(\G{n}{p}).
\end{align*}
\item The same assertion holds for every $H\subseteq\G{n}{p}$ with $\Delta(H)\le(\tau+\delta)\log(np)$.
\end{enumerate}
Conversely, for every $\eta>0$ there is some $\eps_0>0$ such that a.a.s.~there are subgraphs $H_1,H_2\subseteq\G{n}{p}$ satisfying
\begin{align*}
e(H_1)&\le(\tau/2+\eta)n\log(np),&
\Delta(H_2)&\le(\tau+\eta)\log(np),
\end{align*}
and $\chi(\G{n}{p}-H_i)<(1-\eps_0)\gamma\chi(\G{n}{p})$ for $i=1,2$.
\end{theorem}

\begin{proof}
We first state the improper colouring estimates that we need.  We proved in~\cite[Theorem~2]{KaMc10} that, for every integer sequence $t_n$ with $t_n/\log n\to\tau$,
\begin{align*}
\chi^{t_n}(\G{n}{p})\sim\frac{n}{\kappa_p(\tau)\log n}
\end{align*}
a.a.s.  We claim that, for every real sequence $u_n$ with $u_n/\log n\to\tau$,
\begin{align}\label{eqn:global-asymptotic}
\hat{\chi}^{u_n}(\G{n}{p})\sim\frac{n}{\kappa_p(\tau)\log n}
\end{align}
a.a.s.  The upper bound follows from
\begin{align*}
\hat{\chi}^{u_n}(\G{n}{p})\le\chi^{\lfloor u_n\rfloor}(\G{n}{p}),
\end{align*}
since every $\lfloor u_n\rfloor$-improper colouring has at most $n\lfloor u_n\rfloor/2$ monochromatic edges.

For the lower bound, fix $\kappa>\kappa_p(\tau)$ and let $j=\lfloor n/(\kappa\log n)\rfloor$.  Consider a labelled colouring with at most $j$ colours.  By the convexity argument in the proof of Lemma~\ref{lem:lower}, the number of candidate monochromatic pairs is at least $n(\kappa\log n-1)/2$.  Let $N_0=\lfloor n(\kappa\log n-1)/2\rfloor$ and $T=\lfloor nu_n/2\rfloor$.  Since $T/N_0\to\tau/\kappa<p$, we have $T\le pN_0<N_0$ for all sufficiently large $n$.  Applying Lemma~\ref{lem:binomial}, and using $(1-p)^{N_0}=\exp(-N_0\Lambda^*(0))$ for $T=0$, we find that the expected number of globally $u_n$-improper colourings with at most $j$ colours is at most
\begin{align*}
\exp\left(\left[1-\frac{\kappa}{2}
\Lambda^*\left(\frac{\tau}{\kappa}\right)+o(1)\right]n\log n\right).
\end{align*}
This tends to zero exponentially in $n\log n$ by the definition of $\kappa_p(\tau)$.  Letting $\kappa$ decrease to $\kappa_p(\tau)$ proves the lower bound in~\eqref{eqn:global-asymptotic}.

By continuity of $\kappa_p$ and~\eqref{eqn:tau-gamma}, we may choose $\delta>0$ so that both
\begin{align*}
\frac{\kappa_p(0)}{\kappa_p(\tau+\delta)}>(1-\eps/2)\gamma
\quad\text{and}\quad
\frac{\kappa_p(0)}{\kappa_p(\tau+2\delta)}>(1-\eps/2)\gamma.
\end{align*}
Let
\begin{align*}
t_n&=\lfloor(\tau+\delta)\log(np)\rfloor,&
u_n&=(\tau+2\delta)\log(np).
\end{align*}
Since $p$ is fixed, replacing $\log n$ by $\log(np)$ changes only a lower-order term. The two improper colouring estimates above give
\begin{align*}
\frac{\chi^{t_n}(\G{n}{p})}{\chi(\G{n}{p})}
&=(1+o(1))\frac{\kappa_p(0)}{\kappa_p(\tau+\delta)},\\
\frac{\hat\chi^{u_n}(\G{n}{p})}{\chi(\G{n}{p})}
&=(1+o(1))\frac{\kappa_p(0)}{\kappa_p(\tau+2\delta)}.
\end{align*}
Propositions~\ref{prop:localcorr} and~\ref{prop:globalcorr} now give the second and first assertions, respectively; for the first, note that $nu_n/2=(\tau/2+\delta)n\log(np)$.

Finally, $\kappa_p$ is strictly increasing. Hence there is some $\eps_0>0$ such that
\begin{align*}
\frac{\kappa_p(0)}{\kappa_p(\tau+\eta)},\ \frac{\kappa_p(0)}{\kappa_p(\tau+2\eta)}<(1-\eps_0)\gamma.
\end{align*}
Deleting the monochromatic edges of the corresponding improper colourings proves both sharpness assertions.
\end{proof}

The same argument also gives the corresponding assertion for random regular graphs of fixed positive density.

\begin{corollary}\label{cor:sharp,regular,dense}
Fix $p\in(0,1)$ and $\gamma\in(0,1]$, and let $\tau=\tau(p,\gamma)$ be as above.  Suppose $r/n\to p$ and $rn$ is even.  For every $\eps>0$, there is some $\delta>0$ such that the following hold a.a.s.
\begin{enumerate}
\item Every subgraph $H\subseteq\G{n}{r}$ with at most $(\tau/2+\delta)n\log r$ edges satisfies
\begin{align*}
\chi(\G{n}{r}-H)\ge(1-\eps)\gamma\chi(\G{n}{r}).
\end{align*}
\item The same assertion holds for every $H\subseteq\G{n}{r}$ with $\Delta(H)\le(\tau+\delta)\log r$.
\end{enumerate}
Conversely, for every $\eta>0$ there is some $\eps_0>0$ such that a.a.s.~there are subgraphs $H_1,H_2\subseteq\G{n}{r}$ satisfying
\begin{align*}
e(H_1)&\le(\tau/2+\eta)n\log r,&
\Delta(H_2)&\le(\tau+\eta)\log r,
\end{align*}
and $\chi(\G{n}{r}-H_i)<(1-\eps_0)\gamma\chi(\G{n}{r})$ for $i=1,2$.
\end{corollary}

\begin{proof}
We first transfer the two improper colouring estimates used above.  The root $\kappa_p(\tau)$ is jointly continuous in $p$ and $\tau$: on a compact parameter set the unique root of the defining equation can neither approach the boundary $\kappa=\tau/p$ nor tend to infinity, and the assertion follows by taking convergent subsequences.

Fix $s\ge0$ and $\xi>0$.  By this continuity, choose $\delta_0>0$ with $0<p-\delta_0<p+\delta_0<1$ so that
\begin{align*}
\big|\kappa_{p-\delta_0}(s)-\kappa_p(s)\big|,
\ \big|\kappa_{p+\delta_0}(s)-\kappa_p(s)\big|<\xi.
\end{align*}
Choose $0<\eta<1$ so that $\eta p<\delta_0/3$.  Then~\cite[Theorem~1.1]{BeIlMo25}, applied with $\eta$, supplies probabilities $p_-\le p_+$.  Since $r/n\to p$, these probabilities satisfy
\begin{align*}
p-\delta_0<p_-\le p_+<p+\delta_0
\end{align*}
for all sufficiently large $n$.  Extending their sandwich coupling by standard monotone couplings gives
\begin{align*}
\G{n}{p-\delta_0}\subseteq\G{n}{r}\subseteq\G{n}{p+\delta_0}
\end{align*}
a.a.s.  Since $\chi^t$ and $\hat{\chi}^t$ are non-decreasing under the addition of edges, their values for $\G{n}{r}$ lie between the corresponding values for the two binomial graphs.  Let $t_n$ be any integer sequence with $t_n/\log r\to s$, and let $u_n$ be any real sequence with $u_n/\log r\to s$.  Since $\xi>0$ was arbitrary, the fixed-probability estimates and the choice of $\delta_0$ give
\begin{align*}
\chi^{t_n}(\G{n}{r}),\ \hat{\chi}^{u_n}(\G{n}{r})
\sim\frac{n}{\kappa_p(s)\log r}.
\end{align*}
The proof of Theorem~\ref{thm:sharp,dense}, with $\log r$ in place of $\log(np)$, now applies without change.
\end{proof}

For the full range $p\to0$ and $np\to\infty$, a corresponding exact asymptotic formula would follow from matching upper bounds for the local $t$-improper chromatic number, as conjectured in~\cite[Conjecture~4.16]{Kan08} and~\cite{KaMc10}.  These would also give the required global upper bounds, since $\hat{\chi}^t\le\chi^t$.

The asymptotic agreement between the local and global parameters for fixed $p$ need not persist as $p\to1$.  Fix $3/4<a<1$ and write $G=\G{n}{1-n^{-a}}$.  We show that a.a.s.
\begin{align*}
\hat{\chi}^1(G)\le\left(\frac{3}{10}+o(1)\right)n,
\qquad
\chi^1(G)=\left(\frac13+o(1)\right)n.
\end{align*}
Consider $\overline{G}$, which is distributed as $\G{n}{n^{-a}}$.  The expected number of four-cycles in $\overline{G}$ is $O(n^{4-4a})=o(n)$.  Thus the set $B$ of vertices lying in these cycles has order $o(n)$ a.a.s.  Every graph on four vertices with maximum degree at most $1$ has a complement containing a four-cycle, so each class of any $1$-improper colouring contains at most three vertices outside $B$.  Hence $\chi^1(G)\ge(n-|B|)/3$.

For the upper bounds, let $\ell=\lceil n^{(3+a)/4}\rceil=o(n)$.  The probability that some two disjoint sets of $\ell$ vertices have no edge between them in $\overline{G}$ is at most
\begin{align*}
4^n(1-n^{-a})^{\ell^2}\le\exp(2n\log2-n^{-a}\ell^2)=o(1),
\end{align*}
since $n^{-a}\ell^2\ge n^{(3-a)/2}\gg n$.  Thus a.a.s.~$\overline{G}$ has an edge between every two disjoint $\ell$-vertex sets.  For $s\in\{3,4\}$, divide any prescribed vertex set of order $\Theta(n)$ into $s$ equal parts, setting aside fewer than $s$ vertices.  Choose a maximal matching in $\overline{G}$ between each consecutive pair of parts; each matching leaves fewer than $\ell$ vertices unmatched on each side.  All but at most $(s-1)\ell$ vertices of the first part start a path following the matchings through all $s$ parts.  We therefore obtain vertex-disjoint paths on $s$ vertices covering all but $o(n)$ vertices of the chosen set.  These paths need not be induced.

Applying this construction with $s=3$ to $V(G)$ gives triples inducing at most one edge each in $G$.  Using the triples and unused singletons as colour classes yields $\chi^1(G)\le n/3+o(n)$.  For the global bound, choose disjoint sets $V_3,V_4$ of orders $3\lfloor n/5\rfloor$ and $4\lfloor n/10\rfloor$, respectively.  The construction gives $m_3=\lfloor n/5\rfloor-o(n)$ triples in $V_3$ and $m_4=\lfloor n/10\rfloor-o(n)$ quadruples in $V_4$.  Each triple induces at most one edge in $G$, and each quadruple at most three.  We use these triples and quadruples as colour classes, and give each unused vertex a distinct new colour.  This uses $n-2m_3-3m_4=3n/10+o(n)$ colours.  The number of monochromatic edges is at most
\begin{align*}
m_3+3m_4\le\lfloor n/5\rfloor+3\lfloor n/10\rfloor\le n/2,
\end{align*}
as required for global $1$-impropriety.

\section{Concluding remarks}

We have determined the order of chromatic resilience subject to removal of edges in the binomial random graph $\G{n}{p}$ when $np\to\infty$ and $p$ is bounded away from $1$, and in the random regular graph $\G{n}{r}$ when $r\to\infty$ and $r/n$ is bounded away from $1$.  At fixed positive density, we determined in Section~\ref{sec:further} the sharp constants in the bounds on $e(H)$ and $\Delta(H)$ for a prescribed constant-factor reduction in the chromatic number.  Here $H$ is the subgraph whose edges we remove.  The corresponding sharp constants in the sparse regime remain open.

As a by-product of our investigation, we have initiated the probabilistic study of the global $t$-improper chromatic number, a natural parameter that may merit study in other settings for chromatic graph theory.

The bounded-average-degree case, in which $p=d/n$ for fixed $d>0$ or $r$ is fixed, leads to a quantitative refinement of the ordinary colourability problem.  Fix an integer $k\ge3$, and consider $\G{n}{d/n}$.  We ask how many edges we must delete to make this graph $k$-colourable, and how this minimum grows with $d$.  For a graph $G$ on $[n]$, let
\begin{align*}
\operatorname{del}_k(G):=\min\bigl\{e(F):F\subseteq G,\ \chi(G-F)\le k\bigr\}.
\end{align*}
Equivalently, $\operatorname{del}_k(G)$ is the minimum number of monochromatic edges over all $k$-colourings of $G$.  Proposition~\ref{prop:globalcorr} gives
\begin{align*}
\hat{\chi}^t(G)\le k
\quad\Longleftrightarrow\quad
\operatorname{del}_k(G)\le\frac{nt}{2}.
\end{align*}

\begin{samepage}
The following consequence of a result of Bayati, Gamarnik and Tetali~\cite[Theorem~1 and Corollary~1]{BaGaTe13} gives the limiting deletion count.

\begin{proposition}
Fix an integer $k\ge3$.\label{prop:deletion-limit}  There is a deterministic, continuous, non-decreasing and $1/2$-Lipschitz function $\psi_k:(0,\infty)\to[0,\infty)$ such that
\begin{align*}
\operatorname{del}_k(\G{n}{d/n})=(\psi_k(d)+o(1))n
\quad\text{a.a.s.}
\end{align*}
for every fixed $d>0$.  Moreover, $\operatorname{del}_k(\G{n}{d_n/n})/n=\psi_k(d)+o(1)$ a.a.s.~for every sequence $d_n\to d>0$.  Both limits also hold in expectation.

There is a constant $d_k^*>0$ such that $\psi_k(d)=0$ for $0<d\le d_k^*$ and $\psi_k(d)>0$ for $d>d_k^*$.
\end{proposition}
\end{samepage}

\begin{proof}
Bayati, Gamarnik and Tetali~\cite[Theorem~1]{BaGaTe13} prove that, for a uniform random graph $G$ on $[n]$ with $\lfloor dn/2\rfloor$ edges, the maximum number of properly coloured edges in a $k$-colouring, divided by $n$, converges to a deterministic limit.  Subtracting this maximum from the edge count gives the limit for $\operatorname{del}_k(G)/n$ in that model.  The binomial graph is uniform conditional on its edge count, which is $dn/2+o(n)$ a.a.s., so we may couple the two graphs to differ by $o(n)$ edges.  Adding or deleting one edge changes $\operatorname{del}_k$ by at most one, proving the binomial limit.  Monotone coupling at two values of $d$ also gives the stated monotonicity and Lipschitz bound.

The same coupling argument applies to every sequence $d_n\to d>0$.  Both limits also hold in expectation, since $\operatorname{del}_k(G)\le e(G)$ and the second moments of $e(G)/n$ remain bounded.

The assertion concerning $d_k^*$ follows from~\cite[Corollary~1]{BaGaTe13} and continuity of $\psi_k$.
\end{proof}

The value $\psi_k(d)$ is the limiting minimum number of edge deletions per vertex needed for $k$-colourability.  Thus, for fixed $0<d\le d_k^*$, $o(n)$ deletions suffice a.a.s.; for fixed $d>d_k^*$, the minimum is $\Theta(n)$ a.a.s.

\begin{samepage}
\begin{corollary}
Fix an integer $k\ge3$.\label{cor:improper-threshold}  For every fixed $t>0$, there is a unique $b_k(t)>d_k^*$ satisfying $2\psi_k(b_k(t))=t$.  For fixed $0<d<b_k(t)$, we have $\hat{\chi}^t(\G{n}{d/n})\le k$ a.a.s., whereas for fixed $d>b_k(t)$, we have $\hat{\chi}^t(\G{n}{d/n})>k$ a.a.s.  Moreover, $b_k(t)\downarrow d_k^*$ as $t\downarrow0$.
\end{corollary}
\end{samepage}

\begin{proof}
For $0<d_1<d_2$, choose a $k$-colouring of $\G{n}{d_2/n}$ that minimises the number of monochromatic edges.  Retain each edge independently with probability $d_1/d_2$.  The retained graph has distribution $\G{n}{d_1/n}$, and its minimum deletion count is at most the number of surviving monochromatic edges.  Taking expectations and passing to the limit gives
\begin{align*}
\frac{\psi_k(d_1)}{d_1}\le\frac{\psi_k(d_2)}{d_2}.
\end{align*}
Hence $\psi_k$ is strictly increasing on $(d_k^*,\infty)$ and tends to infinity.  Together with continuity and $\psi_k(d_k^*)=0$, this gives the existence and uniqueness of $b_k(t)$ and its stated limit.  Propositions~\ref{prop:globalcorr} and~\ref{prop:deletion-limit} give the two colourability assertions.
\end{proof}

Thus $b_k(t)$ is the limiting average-degree threshold for global $t$-improper $k$-colourability.  The scaling limit gives no conclusion at $d=b_k(t)$.

The function $\psi_k$ combines the questions for different values of $t>0$ into a single quantitative refinement of the Erd\H{o}s--R\'enyi colourability problem~\cite{ErRe60}.  Its relation to ordinary $k$-colourability requires care: $\psi_k(d)=0$ means that $o(n)$ deletions suffice, whereas ordinary $k$-colourability requires no deletions.

Achlioptas and Friedgut~\cite[Theorem~1.1]{AcFr99} proved that ordinary $k$-colourability has a sharp threshold at a sequence of average degrees $d_k(n)$.  Thus, for every fixed $0<\eta<1$,
\begin{align*}
\chi\bigl(\G{n}{(d_k(n)-\eta)/n}\bigr)\le k
\quad\text{and}\quad
\chi\bigl(\G{n}{(d_k(n)+\eta)/n}\bigr)>k
\end{align*}
a.a.s.  We need not assume that $d_k(n)$ converges, although this is conjectured to be the case.  The elementary first-moment bound gives $d_k(n)\le2k\log k+o(1)$, so this sequence is bounded; sharper bounds for large $k$ appear in~\cite[Theorem~1.1]{CoVi16} and~\cite{Coj13}.

Bayati, Gamarnik and Tetali~\cite[Conjecture~1]{BaGaTe13} conjectured exponential decay of the colourability probability above this threshold.  The question that Section~\ref{sec:further} suggests is whether, at every fixed positive distance above the threshold, linearly many edge deletions are necessary for $k$-colourability.  The following proposition relates these assertions to convergence of the threshold to $d_k^*$; we use the simple-graph binomial formulation of the exponential-decay conjecture.

\begin{samepage}
\begin{proposition}
Fix an integer $k\ge3$, and let $d_k(n)$ and $d_k^*$ be as above.\label{prop:threshold-equivalence}  The following statements are equivalent:
\begin{enumerate}
\item\label{item:threshold-convergence}
$d_k(n)\to d_k^*$.
\item\label{item:threshold-resilience}
For every fixed $\eps>0$, there exists $\delta=\delta(k,\eps)>0$ such that
\begin{align*}
\operatorname{del}_k\bigl(\G{n}{(d_k(n)+\eps)/n}\bigr)\ge\delta n
\quad\text{a.a.s.}
\end{align*}
\item\label{item:threshold-decay}
For every fixed $\eps>0$, there exists $\gamma=\gamma(k,\eps)>0$ such that
\begin{align*}
\Pr\Bigl(\chi\bigl(\G{n}{(d_k(n)+\eps)/n}\bigr)\le k\Bigr)\le e^{-\gamma n}
\end{align*}
for all sufficiently large $n$.
\end{enumerate}
\end{proposition}
\end{samepage}

\begin{proof}
For the equivalence of~\ref{item:threshold-resilience} and~\ref{item:threshold-decay}, we transfer the resilience and exponential-decay assertions between the binomial and fixed-edge-count models by monotone coupling and Chernoff bounds, using a smaller positive margin in the source model and incurring an exponentially small coupling error.

Assume~\ref{item:threshold-decay}, and fix $\eps>0$.  Take a uniform random graph $G$ on $[n]$ with $m=\lfloor(d_k(n)+\eps)n/2\rfloor=O(n)$ edges.  Deleting the edges at any fixed $\lfloor\delta n\rfloor$ positions in a uniform random ordering leaves a uniform random graph with $m-\lfloor\delta n\rfloor$ edges.  For fixed $0<\delta<\eps/4$, its average degree exceeds $d_k(n)+\eps/2$ for all sufficiently large $n$.  Exponential decay at a smaller fixed positive margin and a union bound therefore give, for some $\gamma>0$ independent of $\delta$,
\begin{align*}
\Pr\bigl(\operatorname{del}_k(G)\le\lfloor\delta n\rfloor\bigr)
&\le\binom{m}{\lfloor\delta n\rfloor}e^{-\gamma n}\\
&=\exp\bigl((O(\delta\log(1/\delta))-\gamma)n\bigr)=o(1),
\end{align*}
provided $\delta$ is sufficiently small.  For the converse, if $\operatorname{del}_k(G)\ge\delta n$ a.a.s., then its expectation is at least $(\delta+o(1))n$.  The bounded-differences inequality applied to $\operatorname{del}_k(G)$ in the fixed-edge-count model gives an exponentially small bound on $\Pr(\operatorname{del}_k(G)=0)$.  Bayati, Gamarnik and Tetali use this deletion argument and their scaling theorem to deduce threshold convergence from the conjecture.

To deduce~\ref{item:threshold-convergence} from~\ref{item:threshold-resilience}, consider any subsequence along which $d_k(n)\to L$.  The sharp-threshold result gives $\psi_k(d)=0$ for $d<L$, whereas the resilience assertion and Proposition~\ref{prop:deletion-limit} give $\psi_k(d)>0$ for $d>L$.  Hence $L=d_k^*$, and boundedness of $d_k(n)$ proves convergence.  Conversely, if $d_k(n)\to d_k^*$, applying Proposition~\ref{prop:deletion-limit} with $d_n=d_k(n)+\eps$ gives the resilience assertion with any $0<\delta<\psi_k(d_k^*+\eps)$.
\end{proof}

Beyond identifying the transition, we ask how quickly the minimum deletion count grows above it.  This question does not require the exponential-decay conjecture.  The $1/2$-Lipschitz property gives $\psi_k(d_k^*+\eps)\le\eps/2$.  Does a matching linear lower bound hold, or equivalently, is
\begin{align*}
\psi_k(d_k^*+\eps)=\Theta(\eps)
\qquad\text{as }\eps\downarrow0?
\end{align*}

\subsection*{Open access statement}%\mbox{}\\*
For the purpose of open access, a CC BY public copyright license is applied to any Author Accepted Manuscript (AAM) arising from this submission.

\subsection*{Declaration on the use of generative AI}%\mbox{}\\*
We began this project in 2012 and revisited it intermittently thereafter.  The present paper retains the principal questions, results and mathematical ideas of the original manuscript and incorporates subsequent progress on sandwich couplings in its treatment of random regular graphs.

In preparing the present version, the first author used a commercially available, closed-weight generative AI system to improve the exposition, check notation and citations, clarify technical details, and develop and check elementary deductions concerning the concluding remarks and the case $p\to1$.  It was not used to develop the main results or their proof strategies.
%GPT 5.6

\bibliographystyle{abbrv}%abbrv
\bibliography{resilience}

\end{document}